\documentclass[11pt]{amsart}
\newtheorem{theorem}{Theorem}

\newtheorem{lemma}[theorem]{Lemma}
\newtheorem{Def}[theorem]{Definition}

\newcommand{\cc}{\mathbb{C}}
\newcommand{\pp}{\mathbb{P}}
\newcommand{\rr}{\mathbb{R}}
\newcommand{\zz}{\mathbb{Z}}

\title[Rationally convex sets]
{Approximation and the topology of rationally convex sets}

\author{E. S. Zeron}
\address{Depto. Matem\'aticas, CIVESTAV,
Apdo.~Postal 14-740, M\'exico DF, 07000, M\'exico.}
\email{eszeron@math.cinvestav.mx}

\address{Centre de Recherches Math\'ematiques, Universit\'e 
de Montr\'eal, Case pos\-ta\-le 6128, Succursale centre-ville, 
Montr\'eal, H3C 3J7, Canada.}

\date{\today}
\thanks{Research supported by Cinvestav and Conacyt, M\'exico}

\begin{document}
\begin{abstract}
Considering a mapping $g$ holomorphic on a neighbourhood of a rational 
convex set $K\subset\cc^n$, and range into the complex projective space 
$\cc\pp^m$, the main objective of this paper is to show that we can 
uniformly approximate $g$ on $K$ by rational mappings defined from 
$\cc^m$ into $\cc\pp^m$. We only need to ask that the second \v{C}ech 
cohomology group $\check{H}^2(K,\zz)$ vanishes.
\end{abstract}

\maketitle
\section{Introduction}

Let $\cc\pp^m$ be the $m$-complex projective space, composed of all the
complex lines in $\cc^{m+1}$ which pass through the origin. It is well
known that $\cc\pp^m$ is a $m$-complex manifold, and that there exists 
a natural holomorphic projection $\rho_m$ defined from
$\cc^{m+1}\setminus\{0\}$ onto $\cc\pp^m$, which sends 
any point $(z_0,\ldots,z_m)\neq0$ to the complex line
$$\rho_m(z_0,\ldots,z_m)=[z_0,\ldots,z_m]:=\{(z_0t,\ldots,z_mt):t\in\cc\}.$$

In particular, we have that the one-dimensional complex projective 
space $\cc\pp^1$ is the Riemann sphere $\mathcal{S}^2$, and the natural 
holomorphic projection $\rho_1$ is given by $\rho_1(w,z)=[\frac{w}{z},1]$
or $[1,\frac{z}{w}]$. Thus, any rational mapping $p/q$ defined on $\cc^n$ 
may be seen as the composition $\rho_1(p,q)$, where $(p,q)$ is a 
holomorphic polynomial mapping from $\cc^n$ into $\cc^2$. The critical 
set $E$ of $p/q$ is the inverse image $(p,q)^{-1}(0)$, and so $p/q$ is a 
holomorphic mapping defined from $\cc^n\setminus{E}$ into $\mathcal{S}^2$. 
Previous interpretation allows us to extend the notion of rational mapping 
to consider the natural projections $\rho_m$ for $m\geq1$.

\begin{Def}\label{Rat} A rational mapping based on $\cc^n$, and image 
in $\cc\pp^m$, is defined as the composition $\rho_m(P)$, for a given
holomorphic polynomial mapping $P:\cc^n\to\cc^{m+1}$. The critical set 
$E$ of $\rho_m(P)$ is then defined as the inverse image $P^{-1}(0)$, and 
so $\rho_m(P)$ is a holomorphic mapping defined from $\cc^n\setminus{E}$ 
into $\cc\pp^m$.  
\end{Def}

On the other hand, recalling the fundamentals of rational approximation
theory. We have that a compact set $K$ in $\cc^n$ is rationally convex if
for every point $y\in\cc^n\setminus{K}$ there exists a holomorphic
polynomial $p$ such that $p(y)=0$ and $p$ does not vanishes on $K$.  
Besides, it is well known that each function $h:U\to\cc$ holomorphic on a
neighbourhood $U$ of $K$ can be approximated on $K$ by rational functions,
whenever $K$ is rationally convex, see for example \cite{AW} or
\cite{Gam}.  That is, for each $\widehat{\delta}>0$ there exists a
holomorphic rational function $p/q$ such that $K$ does not meet the zero
locus of $q$ and $\bigl|\frac{p(z)}{q(z)}-h(z)\bigr|$ is strictly less
that $\widehat{\delta}$ on $K$. This result automatically drives us to
consider whether the concept of rationally convex sets is strong enough 
as to imply approximation by the kind of rational mappings that we have 
just introduced on Definition~\ref{Rat}. Amazingly, we can get a 
positive answer by adding a simple cohomological condition.

\begin{theorem}{\bf (main theorem)} Let $K$ be a rationally convex set 
in $\cc^n$, and {\rm{Dist}} be a metric on $\cc\pp^m$ which induces 
the topology, with $m,n\geq1$. If the second \v{C}ech cohomology group
$\check{H}^2(K,\zz)$ vanishes. Then, for each $\widehat{\delta}>0$ and 
any mapping $g:U\to\cc\pp^m$ holomorphic on a neighbourhood $U$ of $K$, 
there exists a rational mapping $\rho_m(P)$ defined on $\cc^m$, whose 
critical set does not meet $K$, and such that 
{\rm{Dist}}$\bigl[\rho_m(P(z)),g(z)\bigr]$ 
is less than $\widehat{\delta}$ on $K$.  
\end{theorem}

This result was mainly inspired on the work done by professors Grauert,
Kerner and Oka in \cite{Gra}, \cite{GK} and \cite{Oka}. Besides, in an
early paper \cite{GZ} we have previously analysed the approximation by the
rational mappings described in Definition~\ref{Rat}, we deduced a result
similar to main theorem by using the extra topological condition of being
null-homotopic. We are giving a complete reference to this early result in
section~4.

We shall prove the main theorem in the third section of this paper.  
Moreover, we are going to devote the second section to introduce the
results on cohomology theory that we need for the proof of the main 
theorem. Finally, examples and corollaries are introduced on section~4.

\section{Cohomology}

We strongly recommend the bibliography \cite{AGP}, \cite{Bre} and
\cite{Gray}, for references on homotopy theory; and the bibliography
\cite{Mas}, \cite{Skl} and \cite{Spa}, for references on cohomology
theory.

We consider two main classes of cohomology groups: \v{C}ech and
singular. And these cohomology groups are both isomorphic on smooth
manifolds, and open subsets of $\cc^n$, for these spaces are all locally
contractible, see for example \cite[pp.~166]{Skl} or \cite[pp.~334
and~341]{Spa}. However, there is a very nice example in \cite[p.~77
and~317]{Spa} of a compact set $K\subset\cc$ whose \v{C}ech cohomology
group $\check{H}^1(K,\zz)=\zz$, but its singular cohomology group
$H_s^1(K,\zz)$ vanishes. Now then, given a closed set $E\subset\cc^n$, 
we need \v{C}ech cohomology groups because $\check{H}^*(E,\zz)$ can be
calculated as the direct limit of the sequence $\{\check{H}^*(U,\zz)\}$,
where $U$ runs over a system (directed by inclusions) of open
neighbourhoods of $E$ in $\cc^n$, see for example \cite[chapter~15]{AW},
\cite[p.~348]{Bre}, \cite[p.~145]{Skl} or \cite[p.~327]{Spa}. Whence, the
\v{C}ech cohomology group $\check{H}^*(E,\zz)$ vanishes if and only if for
each element $\xi\in\check{H}^*(U,\zz)$ defined on an open neighbourhood
$U$ of $E$, there exists a second open set $W$ such that
$E\subset{W}\subset{U}$ and the restriction $\xi|_W$ is 
equal to zero in $\check{H}^*(W,\zz)$.

On the other hand, we need singular cohomology groups because of the
following universal result. Let $U\subset\cc^n$ be an open subset which
has the homotopy type of a CW-complex, the singular cohomology group
$H_s^k(U,\zz)$ is then isomorphic to the group of homotopy classes
$[U,Y]$, where $k\geq1$ and $Y$ is an Eilenberg-McLane space of type
$(\zz,k)$, see for example \cite[p.~183]{AGP}, \cite[pp.~488--492]{Bre},
\cite[p.~274]{DP} or \cite[p.~428]{Spa}. Recall that $Y$ is an
Eilenberg-McLane space of type $(\zz,k)$ if every homotopy group
$\pi_*(Y)$ vanishes, but $\pi_k(Y)$ which is equal to $\zz$. Besides,
recall that $[U,Y]$ is the group composed by all the homotopy classes 
of continuous mappings $f:U\to{Y}$.

Combining the ideas presented in previous paragraphs, we may deduce the
following result. Let $E\subset\cc^n$ be a closed set whose \v{C}ech
cohomology group $\check{H}^k(E,\zz)$ vanishes, $k\geq1$. Besides, 
suppose from now on that $E$ has a system (directed by inclusions) of open
neighbourhoods $\{U_{\beta}\}$ in $\cc^n$, where each $U_{\beta}$ has the
homotopy type of a CW-complex and $E=\bigcap_{\beta}U_{\beta}$. Given an
Eilenberg-McLane space $Y$ of type $(\zz,k)$, and since \v{C}ech and
singular cohomology groups are isomorphic on each $U_{\beta}$, we have 
that for every continuous mapping $f:U_{\beta}\to{Y}$ there exists a 
second neighbourhood $U_{\theta}$ such that
$E\subset{U_{\theta}}\subset{U_{\beta}}$ and the restriction
$f|_{U_{\theta}}:U_{\theta}\to{Y}$ is null-homotopic. The main idea behind
this result is to \textit{see} the homotopy class of the mapping $f$ as
an element of $H_s^k(U_{\beta},\zz)\cong\check{H}^k(U_{\beta},\zz)$.

Let us illustrate previous result with a known example. Recalling that
$\cc\setminus\{0\}$ has the homotopy type of the $1$-dimensional sphere
$\mathcal{S}^1$, and that $\mathcal{S}^1$ is a Eilenberg-McLane space of
type $(\zz,1)$, because $\pi_1(\mathcal{S}^1)=\zz$ and
$\pi_k(\mathcal{S}^1)=0$ for every $k\neq1$. We may deduce that for each
non-vanishing continuous function $f:U_{\beta}\to\cc\setminus\{0\}$, and
whenever $\check{H}^1(E,\zz)=0$, there exists a second neighbourhood
$U_{\theta}$ such that the restriction
$f|_{U_{\theta}}:U_{\theta}\to\cc\setminus\{0\}$ has a well defined
continuous logarithm ($f|_{U_{\theta}}$ is null-homotopic). This result
was presented in chapter~15 of \cite{AW} and in section~7, chapter~III, 
of \cite{Gam}.

Coming back to the main theorem of this paper, we are supposing in the
hypotheses that the second \v{C}ech cohomology group $\check{H}^2(K,\zz)$
vanishes, so we need an example of an Eilenberg-McLane space of type
$(\zz,2)$. This example is given by the infinite dimensional complex
projective space $\cc\pp^\infty$, see for example \cite[p.~360]{AGP},
\cite[p.~157]{Gray} or \cite[p.~425]{Spa}. Making again all the
calculations done in previous paragraphs, we may prove the following
lemma.

\begin{lemma} \label{topo} Let $E$ be a closed subset of $\cc^n$ whose
second \v{C}ech cohomology group $\check{H}^2(E,\zz)$ vanishes. Supposing
that $E$ has a system (directed by inclusions) of open neighbourhoods
$\{U_{\beta}\}$ in $\cc^n$, where each $U_{\beta}$ has the homotopy type
of a CW-complex and $E=\bigcap_{\beta}U_{\beta}$. We have that for every
continuous mapping $f:U_{\beta}\to\cc\pp^\infty$ there exists a second
neighbourhood $U_{\theta}$ such that
$E\subset{U_{\theta}}\subset{U_{\beta}}$ and the restriction
$f|_{U_{\theta}}:U_{\theta}\to\cc\pp^\infty$ is null-homotopic.  
\end{lemma}

We may now prove the main theorem of this paper.

\section{Proof of main theorem}

We need to recall some properties about the infinite dimensional complex
projective space $\cc\pp^\infty$. Consider the infinite dimensional space 
$\cc^\infty$ composed of all the complex sequences $(z_0,z_1,\ldots)$, 
where only a finite number of entries $z_k$ is different from zero. This 
space $\cc^\infty$ is naturally endowed with the standard norm 
$\sqrt{\sum}_k|z_k^2|$. The complex projective space $\cc\pp^\infty$ is 
then composed of all the complex lines in $\cc^\infty$ which pass through 
the origin. Besides, there exists a natural projection $\rho_\infty$ 
defined from $\cc^\infty\setminus\{0\}$ onto $\cc\pp^\infty$ which send 
any point $(z_0,z_1,\ldots)\neq0$ to the complex line
$$\rho_\infty(z_0,z_1,\ldots)=[z_0,z_1,\ldots]
:=\{(z_0t,z_1t,\ldots):t\in\cc\}.$$

Finally, it is easy to calculate that $\rho_\infty$ induces a locally
trivial fibre bundle in $\cc^\infty\setminus\{0\}$, with base on
$\cc\pp^\infty$ and fibre $\cc\setminus\{0\}$, see \cite[p.~360]{AGP}.  
We may cover $\cc\pp^\infty$ with open sets $W_k$ composed of all points
in $\cc\pp^\infty$ whose $k$-entry is equal to one. The open set $W_0$ 
is equal to $\{[1,y]:y\in\cc^\infty\}$, for example. It is now easy to
calculate that $\rho_\infty$ induces a trivial fibre bundle on each
$\rho^{-1}_{\infty}(W_k)$, with base on $W_k$ and fibre
$\cc\setminus\{0\}$. Actually, for every $m\geq1$, we have that 
the projection $\rho_m$ induces a locally trivial fibre bundle 
in $\cc^{m+1}\setminus\{0\}$, with base on $\cc\pp^m$ and fibre
$\cc\setminus\{0\}$, as well. We may now prove the main theorem.

\begin{proof}{\bf (main theorem).} Define the open 
rational polyhedra $V_\beta$ in $\cc^n$ by the formula,
\begin{equation}\label{eq0}
V_\beta\,:=\,\{z\in\cc^n:\,|p_j(z)|<1,\forall{j},
\hbox{~and~}|q_k(z)|>1,\forall{k}\},
\end{equation} 
for some given collections of holomorphic polynomials $\{p_j\}$ and
$\{q_k\}$ in $\cc^n$. We have that each rational polyhedron $V_\beta$ is 
an open Stein subset of $\cc^n$, so $V_\beta$ has the homotopy type of 
a CW-complex \cite[p.~39]{Mil}. Moreover, given a rationally convex set
$K\subset\cc^n$, it is easy to see that the family of all open rational
polyhedra $V_\beta$ which contain to $K$ form a system of neighbourhoods
in $\cc^n$, and that $K$ is equal to the intersection
$\bigcap_{K\subset{V}_\beta}V_\beta$.

On the other hand, let $g:U\to\cc\pp^m$ be any mapping holomorphic on a
neighbourhood $U$ of $K$. We can obviously extend this mapping to a second
one with range on $\cc\pp^\infty$, we only need to set the first $m+1$
entries equal to the entries of $g=[g_0,\ldots,g_m]$ and the rest of 
them equal to zero. That is, define $\breve{g}:U\to\cc\pp^\infty$ by
\begin{equation}\label{eq1}
\breve{g}(z)\,:=\,[g_0(z),\ldots,g_m(z),0,\ldots,0,\ldots]. 
\end{equation}

We may find an open rational polyhedron $V_\beta$ such that
$K\subset{V_\beta}\subset{U}$. Besides, recalling that the \v{C}ech
cohomology group $\check{H}^2(K,\zz)$ vanishes because the given
hypotheses, and considering Lemma~\ref{topo}, we may even find a 
rational polyhedron $V_\beta$ such that the restriction
$\breve{g}|_{V_\beta}:V_\beta\to\cc\pp^\infty$ is null-homotopic. Let
$I=[0,1]$ the unit closed interval in the real line. There exists then 
a continuous mapping $G$ from $V_\beta\times{I}$ into $\cc\pp^\infty$ 
such that $G(z,1)=\breve{g}(z)$ and $G(z,0)=[1,0,\ldots]$, for every
$z\in{V_\beta}$. We have the following commutative diagram,
$$\begin{array}{rcl}
V_\beta&\stackrel{c}{\longrightarrow}&\cc^\infty\setminus\{0\}\\
\downarrow^{\jmath}\;&&\;\downarrow^{\rho_\infty}\\
V_\beta\times{I}&\stackrel{G}{\longrightarrow}&\cc\pp^\infty,\\
\end{array}$$ 
where $c(z)=(1,0,\dots)$ is a constant mapping and $\jmath(z)=(z,0)$ is
the natural inclusion. We know that the projection $\rho_\infty$ induces 
a locally trivial fibre bundle on $\cc^\infty\setminus\{0\}$, with base
$\cc\pp^\infty$ and fibre $\cc\setminus\{0\}$. This fibre bundle has the
homotopy lifting property, see for example \cite[pp.~62 and~67]{DP},
\cite[p.~87]{Gray} or \cite[p.~96]{Spa}. Hence, there exists a continuous
mapping $F$ from $V_\beta\times{I}$ into $\cc^\infty\setminus\{0\}$ such
that $\rho_\infty(F)$ is identically equal to $G$ on $V_\beta\times{I}$.
Recalling equation~(\ref{eq1}), where $\breve{g}(z)=G(z,1)$ was defined,
we can deduce that $F(z,1)$ has the form,
$$F(z,1):=\bigl(F_0(z,1),\ldots,F_m(z,1),0,\ldots,0,\ldots\bigr).$$

We may then introduce a new continuous mapping $f$ defined from $V_\beta$
into $\cc^{m+1}\setminus\{0\}$, by removing the last entries of $F(z,1)$
equal to zero; that is,
$$f(z)\,:=\,\bigl(F_0(z,1),\ldots,F_m(z,1)\bigr)
\quad\hbox{for every}\quad z\in V_\beta.$$

It is easy to see that $\rho_m(f(z))=g(z)$ for every $z\in{V_\beta}$.  
The main objective of previous calculations was the construction of the
continuous mapping $f$ described above. Actually, we could have showed the
existence of such a mapping $f$ by using the results on obstruction theory
described in \cite[p.~507]{Bre} and \cite[p.~447]{Spa}. However, we think
that the procedure followed in previous paragraphs is much more simple and
illustrative.

Nevertheless, we look for a holomorphic (not only continuous) 
mapping $h$ from $V_\beta$ into $\cc^{m+1}\setminus\{0\}$ such that
$\rho_m(h(z))=g(z)$ for every $z\in{V_\beta}$. We shall construct this
holomorphic mapping $h$ by using Oka's results on the second Cousin
problem. Define the space,
$$M\,:=\,\{(z,w)\in{V_\beta}\times\cc^{m+1}:\,g(z)=\rho_m(w),\,w\neq0\}.$$

It is easy to deduce that $M$ is an analytic space because $g$ and
$\rho_m$ are both holomorphic mappings. Moreover, we also have the
following commutative diagram,
$$\begin{array}{lcl}
M&\stackrel{\eta_2}{\longrightarrow}&\cc^{m+1}\setminus\{0\}\\
\downarrow^{\eta_1}&&\;\downarrow^{\rho_m}\\
V_\beta&\stackrel{g}{\longrightarrow}&\cc\pp^m,\\ 
\end{array}$$ 
where $\eta_1(z,w)=z$ and $\eta_2(z,w)=w$ are the basic projections. It 
is easy to prove that $\eta_1$ induces a locally trivial fibre bundle 
in $M$, with Stein base on $V_\beta$ and Stein fibre $\cc\setminus\{0\}$.  
This fibre bundle $M\stackrel{\eta_1}{\longrightarrow}V_\beta$ is the
\textit{pull back} of the fibre bundle induced by $\rho_m$ on
$\cc^{m+1}\setminus\{0\}$. Now then, recalling the continuous mapping 
$f$ defined above, we automatically have that $z\mapsto(z,f(z))$ is a
continuous section of the fibre bundle
$M\stackrel{\eta_1}{\longrightarrow}V_\beta$, because $g(z)=\rho_m(f(z))$
for every $z\in{V_\beta}$. Oka's results on the second Cousin problem
imply that $z\mapsto(z,f(z))$ is homotopic to a holomorphic section
$z\to(z,h(z))$, because $V_\beta$ is Stein, see for example \cite{For},
\cite{FP}, \cite{GK} or \cite{Oka}. Hence, there exists a holomorphic
mapping $h:V_\beta\to\cc^{m+1}\setminus\{0\}$ such that $g(z)$ is equal 
to $\rho_m(h(z))$ for every $z\in{V_\beta}$.

On the other hand, Let $W$ be an open subset of $\cc^{m+1}$ which contains
the compact image $h(K)$. Besides, suppose that the closure $\overline{W}$
is compact and does not contains the origin. Notice that the projection
$\rho_m$ from $\cc^{m+1}\setminus\{0\}$ into $\cc\pp^m$ is continuous with
respect to the metric Dist, which induces the topology, so $\rho_m$ is
also uniformly continuous on $\overline{W}$. Express $h=(h_0,\ldots,h_m)$
as a vector in $\cc^{m+1}$. There are $m+1$ small enough constants
$\delta_k>0$ such that, given $z\in{K}$ and $w\in\cc^{m+1}$,
\begin{eqnarray}\label{eq2}
\hbox{Dist}\bigl[\rho_m(h(z)),\rho_m(w)\bigr]<\widehat{\delta}
&\hbox{and}&w\in\overline{W},\\
\nonumber\hbox{whenever}\quad|h_k(z)-w_k|<\delta_k
&\hbox{for}&0\leq{k}\leq{m}.
\end{eqnarray}

Recalling that $K$ is rationally convex, we may find $m+1$ rational
functions $w_k=\frac{p_k}{q_k}$ defined on $\cc^n$, such that $K$ meets 
the zero locus of no $q_k$, for every $0\leq{k}\leq{m}$, and the absolute 
value $\bigl|h_k(z)-\frac{p_k}{q_k}(z)\bigr|$ is strictly less than 
$\delta_k$ on $K$. Consider the polynomial mapping $P:\cc^n\to\cc^{m+1}$ 
given by
$$\textstyle{P\,:=\,\left(\frac{p_0}{q_0}\prod_{k=0}^mq_k,
\ldots,\frac{p_m}{q_m}\prod_{k=0}^mq_k\right)}.$$

It is easy to deduce that neither
$\bigl(\frac{p_0}{q_0},\ldots,\frac{p_m}{q_m}\bigr)$ nor the product
$\prod_{k=0}^mq_k$ can vanishes on $K$, because $0\not\in\overline{W}$ and
equation~\ref{eq2}. Thus, the compact set $K$ does not meet the critical
set $P^{-1}(0)$ of the rational mapping $\rho_m(P)$. We may also deduce
that $\rho_m(P)$ is equal to 
$\rho_m\bigl(\frac{p_0}{q_0},\ldots,\frac{p_m}{q_m}\bigr)$ on $K$. 
Therefore, recalling equation~\ref{eq2}, and that $\rho_m(h(z))$ is equal 
to $g(z)$ for every $z\in{K}$, we get the result that we look for: the 
distance Dist$\bigl[g(z),\rho_m(P(z))\bigr]$ is strictly less that 
$\widehat{\delta}$ on $K$. 
\end{proof}

\section{Examples and applications}

We shall conclude this paper by noticing that the cohomological condition
$\check{H}^2(K,\zz)=0$ in the hypotheses of the main theorem is not a
trivial condition. Consider the standard two-sphere in $\rr^3$,
$$\mathcal{S}^2\,:=\,\{(a,b,c)\in\rr^3:\,a^2+b^2+c^2=1\}.$$

We can obviously analyse $\mathcal{S}^2$ as a subset of $\cc^3$, 
embedding it into the real axis; and it is easy to see that 
$\mathcal{S}^2$ is rationally convex in $\cc^3$. Besides, the groups
$\check{H}^2(\mathcal{S}^2,\zz)$ and $H_s^2(\mathcal{S}^2,\zz)$ are both
isomorphic to $\zz$. Finally, consider the following open neighbourhood
$U$ of $\mathcal{S}^2$,
$$U\,:=\,\{(x,y,z)\in\cc^3:\,-\pi<\arg(x^2+y^2+z^2)<\pi\}.$$

And the holomorphic mapping $g:U\to\cc\pp^1$, where $\sqrt{1}=1$,
$$\textstyle{g(x,y,z)=
\bigl[\frac{x+iy}{\surd(x^2+y^2+z^2)+z},1\bigr]\quad\hbox{or}
\quad\bigl[1,\frac{x-iy}{\surd(x^2+y^2+z^2)-z}\bigr].}$$

It is easy to see that the restriction $g|_{\mathcal{S}^2}$ is 
the identity mapping from $\mathcal{S}^2$ onto $\cc\pp^1$, because
$g(a,b,c)=\bigl[\frac{a+ib}{1+c},1\bigr]$ is the stereographic projection,
so $g|_{\mathcal{S}^2}$ is not null-homotopic. Moreover, we assert that
$g$ cannot be approximated on $\mathcal{S}^2$ by rational mappings
$\rho_1(P)$. That is, there exists a fix constant $\beta_g>0$, such that 
for every rational mapping $\rho_1(P)$ holomorphic on $\mathcal{S}^2$ 
there is a point $\breve{w}\in\mathcal{S}^2$ with
$\hbox{Dist}\bigl[\rho_1(P(\breve{w})),g(\breve{w})\bigr]$ 
greater than $\beta_g$.

Let $P:\cc^3\to\cc^2$ be any polynomial mapping whose fibre $P^{-1}(0)$
does not meet $\mathcal{S}^2$. We can deduce that the restriction
$P|_{\mathcal{S}^2}$ defined from $\mathcal{S}^2$ into
$\cc^2\setminus\{0\}$ is null-homotopic, because $\cc^2\setminus\{0\}$ 
has the homotopy type of the three-sphere $\mathcal{S}^3$ and the second
homotopy group $\pi_2(\mathcal{S}^3)$ vanishes. Whence, we also have that
the restriction $\rho_1(P)|_{\mathcal{S}^2}$ is null-homotopic for every
rational mapping $\rho_1(P)$ defined on $\cc^3$, and whose critical set
$P^{-1}(0)$ does not meet $\mathcal{S}^2$. Finally, since
$g|_{\mathcal{S}^2}$ is not null-homotopic, we can conclude that there
exists a fix constant $\beta_g>0$, such that for every rational mapping
$\rho_1(P)$ holomorphic on $\mathcal{S}^2$ there is a point 
$\breve{w}\in\mathcal{S}^2$ with
$\hbox{Dist}\bigl[\rho_1(P(\breve{w})),g(\breve{w})\bigr]$ greater 
than $\beta_g$. We only need to recall that $\cc\pp^1$ is an absolute 
neighbourhood retract \cite[pp.~332 and~339]{Kur}, and to apply the 
following lemma which was originally presented in \cite{GZ}.

\begin{lemma}\label{Lm1} Let $X$ and $(Y,\hbox{\rm{d}})$ be two metric
spaces, such that $X$ is compact and $Y$ is an absolute neighbourhood
retract. Then, given a fixed continuous mapping $g:X\to{Y}$, there exists
a constant $\beta_g>0$ such that: every continuous mapping $f:X\to{Y}$ is
homotopic to $g$, whenever {\rm{d}}$[g(x),f(x)]$ is less than $\beta_g$ 
for every $x\in{X}$. 
\end{lemma}

\begin{proof} Let $Y^X$ be the topological space composed of all the
continuous mappings $f:X\to{Y}$, and endowed with the compact-open
topology. Since $X$ is compact and $(Y,\hbox{d})$ is metric, the
compact-open topology of $Y^X$ is induced by the metric
$$\hbox{D}[f_1,f_2]:=\sup\big\{\hbox{d}[f_1(x),f_2(x)]:\,x\in{X}\bigr\},$$
for any two mappings $f_1$ and $f_2$ in $Y^X$ (see for example
\cite[p.~89]{Kur}). The space $Y^X$ is locally arcwise connected and an
absolute neighbourhood retract, because $Y$ is an absolute neighbourhood
retract (see \cite[pp.~339--340]{Kur}). Whence, there exists a fixed
constant $\beta_g>0$ such that the open ball in $Y^X$ with centre in $g$
and radius $\beta_g$ is contained in an arcwise connected neighbourhood 
of $g$. That is, for every continuous mapping $f:X\to{Y}$ with
D$[f,g]<\beta_g$, there exists an arc in $Y^X$ whose end points 
are $f$ and $g$; and so the mappings $g$ and $f$ are homotopic.  
\end{proof}

Finally, the ideas introduced in the proof of the main theorem may be used
to show several versions of this theorem. For example, given a closed set
$E$ in $\cc^n$, a continuous function $f:E\to\cc$ can be tangentially 
approximated by meromorphic ones if for every strictly positive continuous 
function $\hat{\epsilon}:E\to\rr$ there exists a pair of holomorphic 
functions $\phi$ and $\psi$ defined from $\cc^n$ into $\cc$, such 
that $E$ does not meet the zero locus of $\psi$ and
$\bigl|f(z)-\frac{\phi(z)}{\psi(z)}\bigr|$ is less 
than $\hat{\epsilon}(z)$ for every $z\in{E}$. We 
may show the following result now.

\begin{theorem} Let $E$ be a closed set in $\cc^n$ whose second \v{C}ech
cohomology group $\check{H}^2(E,\zz)$ vanishes, and such that every
continuous function $f:E\to\cc$ can be tangentially approximated by
meromorphic functions. Suppose that {\rm{Dist}} is a metric on 
$\cc\pp^m$ which induces the topology, and that $E$ has a system of open 
neighbourhoods $\{V_{\beta}\}$ in $\cc^n$, where each $V_{\beta}$ has the
homotopy type of a CW-complex and $E=\bigcap_{\beta}V_{\beta}$.

For every pair of continuous mappings $\xi:E\to\cc\pp^m$ and 
$\hat{\epsilon}:E\to\rr$, with $\hat{\epsilon}(z)>0$, there exists 
a holomorphic mapping $H$ defined from $\cc^n$ into $\cc^{m+1}$, 
such that $E$ does not meet the zero locus of $H$ and 
{\rm{Dist}}$\bigl[\rho_m(H(z)),\xi(z)\bigr]$ is less 
than $\hat{\epsilon}(z)$ for every $z\in{E}$. 
\end{theorem}

\begin{proof} We only give a sketch of this proof, for it is essentially
the same one presented in section three. Firstly, we have that $\cc\pp^m$
is an absolute neighbourhood retract, for it is homeomorphic to a compact
polyhedron, see \cite[pp.~332 and~339]{Kur}. Therefore, there exists a
continuous mapping $g:U\to\cc\pp^m$ defined on an open neighbourhood $U$
of $E$ and such that $g(z)=\xi(z)$ for every $z\in{E}$. Besides, following
equation~(\ref{eq1}), we can extend $g$ to a continuous mapping
$\breve{g}$ defined from $U$ into $\cc\pp^\infty$.

Considering Lemma~\ref{topo}, there exists an open neighbourhood 
$V_\beta$ such that $E\subset{V_\beta}\subset{U}$ and the restriction 
$\breve{g}|_{V_\beta}:V_\beta\to\cc\pp^\infty$ is null-homotopic. 
Following the ideas presented in pages~4 and~5 of this paper, we can 
build a second continuous (not necessarily holomorphic) mapping $h$ from 
$V_\beta$ into $\cc^{m+1}\setminus\{0\}$ such that $\rho_m(h(z))$ is equal 
to $g(z)=\xi(z)$ for every $z\in{E}$. Notice that $\rho_m$ is continuous 
with respect to the metric Dist, which induces the topology. Express 
$h=(h_0,\ldots,h_m)$ as a vector in $\cc^{m+1}\setminus\{0\}$. There are 
$m+1$ strictly positive continuous functions $\delta_k:E\to\rr$ such that, 
given $z\in{E}$ and $w\in\cc^{m+1}$,
\begin{eqnarray}\label{eq7}
\hbox{Dist}\left[\rho_m(h(z)),\rho_m(w)\right]<\hat{\epsilon}(z)
\quad\hbox{whenever}&&\\
\nonumber|h_k(z)-w_k|<\delta_k(z)\quad\hbox{for}\quad0\leq{k}\leq{m}.
&&\end{eqnarray} 

Recalling that every continuous function can be tangentially approximated
by meromorphic functions on $E$, we may find holomorphic functions
$\phi_k$ and $\psi_k$ defined from $\cc^n$ into $\cc$ such that $E$ meets
the zero locus of no $\psi_k$, for $0\leq{k}\leq{m}$, and the absolute
value $\bigl|h_k(z)-\frac{\phi_k}{\psi_k}(z)\bigr|$ is strictly less 
than $\delta_k(z)$ for every $z\in{E}$. Moreover, we can even choose the
functions $\delta_k$ in such a way that each $\delta_k(z)$ is strictly
less than $\max\bigl\{\frac{|h_j(z)|}{2}\bigr\}$ for all $z\in{E}$. 
Hence, we have that neither
$\bigl(\frac{\phi_0}{\psi_0},\ldots,\frac{\phi_m}{\psi_m}\bigr)$ 
nor the product $\prod_{k=0}^m\psi_k$ can vanishes on $E$; and so, 
the set $E$ does not meet the zero locus of the holomorphic mapping
$H:\cc^n\to\cc^{m+1}$ defined by,
$$\textstyle{H\,:=\,\left(\frac{\phi_0}{\psi_0}\prod_{k=0}^m\psi_k,
\ldots,\frac{\phi_m}{\psi_m}\prod_{k=0}^m\psi_k\right).}$$

It is easy to deduce that $\rho_m(H)$ is equal to
$\rho_m\bigl(\frac{\phi_0}{\psi_0},\ldots,\frac{\phi_m}{\psi_m}\bigr)$ 
on $E$. Whence, recalling equation~\ref{eq7}, and that $\rho_m(h(z))$ 
is equal to $g(z)=\xi(z)$ for every $z\in{E}$, we get the result that 
we look for: Dist$\left[\xi(z),\rho_m(H(z))\right]$ is strictly less 
than $\hat{\epsilon}(z)$ on $E$. 
\end{proof}

We have already proved a weaker version of this theorem in 
\cite{GZ}. We used there the stronger hypotheses that $E$ is 
compact and $\xi:E\to\cc\pp^m$ is null-homotopic, instead 
of asking $\check{H}^2(E,\zz)=0$.

\bibliographystyle{amsplain}

\begin{thebibliography}{10}

\bibitem{AGP} M.A. Aguilar, S. Gitler and C. Prieto. 
\textit{Topolog{\'\i}a algebraica, un enfoque homot\'opico}.
McGraw-Hill, M\'exico, 1998.

\bibitem{AW} H. Alexander and J. Wermer. \textit{Several Complex
Variables and Banach Algebras, third edition}, (Graduate Texts
in Mathematics, 35). Springer-Verlag, New York, 1998.

\bibitem{Bre} Glen E. Bredon. \textit{Topology and geometry},
(Graduate Texts in Mathematics, 139). Springer-Verlag, 
New York, 1993.

\bibitem{DP} C.T.J. Dodson and P.E. Parker. \textit{A user's
guide to algebraic topology}. Kluwer Academic Publishers, 
Dordrecht the Netherlands, 1997.

\bibitem{For} F. Forstneri\v{c}. Oka's principle for sections 
of subelliptic submersions. \textit{Math. Z.} \textbf{241} 
(2002), pp.~527--551.

\bibitem{FP} F. Forstneri\v{c} and J. Prezelj. Oka's principle 
for holomorphic submersions with sprays. \textit{Math. Ann.} 
\textbf{322} (2002), pp.~633--666. 

\bibitem{Gam} T.W. Gamelin. \textit{Uniform algebras}.
Prentice-Hall, Englewood Cliffs N.J., 1969.

\bibitem{GZ} Paul M. Gauthier and E. S. Zeron. Approximation by Rational
Mappings, via Homotopy Theory. \textit{Canad. Math. Bull.} in press.

\bibitem{Gra} H. Grauert. Holomorphe Funktionen mit Werten in komplexen 
Lieschen Gruppen. \textit{Math. Ann.} \textbf{133} (1957), pp.~450--472.

\bibitem{GK} H. Grauert and H. Kerner. Approximation von holomorphen
Schnittfl\"achen in Faserb\"undeln mit homogener Faser.
\textit{Arch. Math.} \textbf{14} (1963), pp.~328--333.

\bibitem{Gray} B. Gray. \textit{Homotopy theory, an introduction
to algebraic topology}. Academic Press, New York, 1975.

\bibitem{HW} L. H\"ormander and J. Wermer. Uniform approximation
on compact sets in $\cc^n$. \textit{Math. Scand.} \textbf{23}
(1968), pp.~5--21.

\bibitem{Kur} K. Kuratowski. \textit{Topology, Vol. II}.
Academic Press, New York and London, 1968.

\bibitem{Mas} W.S. Massey. \textit{Homology and cohomology theory}.
Marcel Dekker, New York, 1978.

\bibitem{Mil} J. Milnor. \textit{Morse theory}, (Annals of Mathematics
Studies, 51). Princeton University Press, Princeton N.J., 1963.

\bibitem{NW} R. Nirenberg and R. O. Wells Jr.. Approximation
theorems on differentiable submanifolds of a complex manifold.
{\it Trans. Amer. Math. Soc.} \textbf{142} (1969) pp.~15--35.

\bibitem{Oka} K. Oka. Sur les fonctions des plusieurs variables III:
Deuxi\`eme probl\`eme de Cousin. \textit{J. Sc. Hiroshima Univ.} 
\textbf{9} (1939), pp.~7--19.

\bibitem{Skl} E. G. Sklyarenko. {\it Homology and cohomology theories
of general spaces}. (General topology, II, Encyclopaedia Math. Sci.,
50), Springer, Berlin, 1996, pp.~119--256.

\bibitem{Spa} E. H. Spanier. {\it Algebraic topology}. McGraw-Hill,
New York-Toronto, 1966.

\end{thebibliography}

\end{document}